\documentclass[11pt]{amsart}
\usepackage{etoolbox}
\usepackage{aliascnt,amsmath,amssymb,amsthm,mathtools,xcolor}

\usepackage[shortlabels]{enumitem}
\usepackage[alphabetic,abbrev]{amsrefs}
\usepackage[margin=30truemm]{geometry}
\definecolor{mylinkcolor}{RGB}{16,156,81}
\definecolor{mycitecolor}{RGB}{20,80,140}
\usepackage{hyperref}
\usepackage{nicefrac}
\usepackage[nameinlink]{cleveref}
\hypersetup{colorlinks=true,linkcolor=mylinkcolor,citecolor=mycitecolor,urlcolor=mycitecolor}

\newtheorem{thm}{Theorem}
\newaliascnt{prop}{thm}
\newtheorem{prop}[prop]{Proposition}
\aliascntresetthe{prop}
\newaliascnt{lem}{thm}
\newtheorem{lem}[lem]{Lemma}
\aliascntresetthe{lem}
\theoremstyle{remark}
\newaliascnt{rk}{thm}
\newtheorem{rk}[rk]{Remark}
\aliascntresetthe{rk}
\crefname{thm}{Theorem}{Theorems}
\crefname{prop}{Proposition}{Propositions}
\crefname{lem}{Lemma}{Lemmas}
\crefname{rk}{Remark}{Remarks}
\crefname{table}{Table}{Tables}
\crefname{section}{Section}{Sections}

\newcommand{\Q}{\mathbb Q}
\newcommand{\Z}{\mathbb Z}
\newcommand{\C}{\mathbb C}
\newcommand{\F}{\mathbb F_2}
\newcommand{\Zt}{\Z_{(2)}}
\newcommand{\Gf}[1]{G^{f}_{#1}}
\newcommand{\fmes}{\mathcal G^{\!f}}

\newcommand{\cT}{\mathcal T}
\DeclareMathOperator{\spanQ}{span_{\Q}}
\DeclareMathOperator{\spanF}{span_{\F}}
\def\={\: =\: }
\def\+{\,+\,}
\def\meno{\,-\,}
\def\defis{\: :=\: }

\title[Products of at most two Eisenstein series]{An elementary proof that the space of modular forms is spanned by products of at most two Eisenstein series}
\author{Henrik Bachmann}
\address{Graduate School of Mathematics, Nagoya University, Nagoya, Japan.}
\email{henrik.bachmann@math.nagoya-u.ac.jp}
\author{Jan-Willem M. van Ittersum}
\address{Korteweg--de Vries Institute for Mathematics, University of Amsterdam, Amsterdam, The Netherlands.}
\email{j.w.m.vanittersum@uva.nl}

\date{\today}

\makeatletter
\patchcmd{\ps@firstpage}{\hfil\thepage\hfil}{\hfil\raisebox{-40pt}[0pt][0pt]{\thepage}\hfil}{}{\PackageError{main}{patch of firstpage failed}{}}
\makeatother

\begin{document}

\begin{abstract}
We give an elementary proof that every product of three Eisenstein series of level one is a linear combination of products of two, so that in every weight the space of modular forms is spanned by products of at most two Eisenstein series. The proof uses only the quadratic recurrence for the Eisenstein series, and hence also applies to formal Eisenstein series. Its main ingredient is that the Laurent coefficients of the Weierstrass $\wp$-function, suitably normalized, reduce modulo~$2$ to the power sums of the roots of the associated cubic.
\end{abstract}
\dedicatory{To Don Zagier, on the occasion of his 75th birthday}

\maketitle
\enlargethispage{40pt}

\section{Introduction}
For even $k\geq4\mspace{1mu}$, we consider the Eisenstein series
\[
G_k(\tau)\defis\frac{1}{2\,(2\pi i)^k}\sum_{(m,n)\neq(0,0)}\frac{1}{(m\tau+n)^k}\=-\frac{B_k}{2\,k!}\+\frac{1}{(k-1)!}\sum_{m,n\geq1}m^{k-1}q^{mn}\quad (q=e^{2\pi i\tau}),
\]
where $B_k$ is the $k$th Bernoulli number. We view $G_k$ as an element of $\Q[\![q]\!]\mspace{1mu}$. The Eisenstein series satisfy the quadratic relations
\begin{equation}\label{eq:recurrence-G}
G_k\=\frac{12}{(k-6)(k-1)(k+1)}\sum_{\substack{a+b=k\\a,b\geq4\text{ even}}}(a-1)(b-1)\,G_a\,G_b\qquad(k\geq8)\,,
\end{equation}
the first three of which are $G_8=\frac67G_4^2\mspace{1mu}$, $G_{10}=\frac{10}{11}G_4G_6$ and $G_{12}=\frac{72}{143}G_4^3+\frac{50}{143}G_6^2\mspace{1mu}$. Such identities are usually proved by comparing modular forms in a space of known dimension. Skoruppa gave a short, elementary proof of them in~\cite{S}. In this note, we show the following result in an elementary way and by just using \eqref{eq:recurrence-G}.

\begin{thm}\label{thm:triple}
For even $a,b,c\geq4\mspace{1mu}$, we have
\[
G_a\,G_b\,G_c\;\in\;\spanQ\{G_r\,G_s\::\:r,s\geq4\text{ even},\ r+s=a+b+c\}\,.
\]
\end{thm}

For example,
\[
G_4^3\=\frac76\,G_4\,G_8\,,\qquad G_6^3\=\frac{143}{50}\,G_6\,G_{12}\meno\frac{231}{125}\,G_8\,G_{10}\,.
\]
\Cref{thm:triple} is well known. Since the algebra of modular forms of level one with rational Fourier coefficients is $\Q[G_4,G_6]\mspace{1mu}$, it is equivalent, by induction on the number of factors, to the statement that for even $k\geq8$ the space~$M_k$ of such modular forms of weight~$k$ satisfies
\begin{equation}\label{eq:Mk}
M_k\=\spanQ\{G_a\,G_{k-a}\::\:4\leq a\leq k-4\text{ even}\}\,.
\end{equation}
(The Eisenstein series $G_k$ lies in the right-hand side by~\eqref{eq:recurrence-G}.) The classical proof of~\eqref{eq:Mk} is analytic: by Rankin's method, the Petersson product of a normalized Hecke eigenform~$f$ with~$G_aG_{k-a}$ is, up to a nonzero factor independent of~$f$, the product of the $L$-values $L(f,k-1)$ and $L(f,a)$. Hence, if a cusp form $g$ is orthogonal to all the products $G_aG_{k-a}\mspace{1mu}$, multiplying each of its components in the basis of Hecke eigenforms~$f$ by $L(f,k-1)$ gives a cusp form whose odd periods all vanish, and which is therefore zero by the Eichler--Shimura theorem (see~\cite{R}, \cite{KZ}*{Sections~1.1 and~1.4} and the sketch in~\cite{S}). 

In~\cite{BI}*{before Corollary~5.9} we ask for a combinatorial proof of the analogue of~\eqref{eq:Mk} for the formal Eisenstein series~$\Gf{k}\mspace{1mu}$. Formal Eisenstein series are elements of an algebra~$\fmes$ of words in formal symbols in which each word is identified with its image under a certain involution. They satisfy~\eqref{eq:recurrence-G} 
and $\Gf{k}\mapsto G_k$ under an algebra homomorphism $\fmes\to\Q[\![q]\!]\mspace{1mu}$. Since the $\Gf{k}$ are not functions, analytic tools, like Rankin's method, are not available in~$\fmes\mspace{1mu}$. However, as our proof uses only~\eqref{eq:recurrence-G}, it follows that the analogues of \cref{thm:triple} and~\eqref{eq:Mk} for the formal Eisenstein series also hold.

Our result also makes an elementary description of all linear relations among the~$G_k$ and the products $G_aG_{k-a}$ possible. Hirose, Sato and Tasaka \cite{HST}*{Theorem~1} proved, by partial fraction decomposition and without using modularity, a family of such relations, indexed by triples of positive integers. To show that this family gives all linear relations, and to obtain an explicit basis of~$M_k\mspace{1mu}$, they used~\eqref{eq:Mk} (see also \cite{GKZ}*{Theorem~5}). 

Analogues of~\eqref{eq:Mk} for congruence subgroups were proved in~\cites{BG2,KoM,DN,RX}, again by Rankin's method and the nonvanishing of certain $L$-values. It would be interesting to find elementary proofs in higher level as well.

\subsection*{A statement about polynomials}
To prove \cref{thm:triple}, we replace $G_4$ and~$G_6$ by variables. Thus, we give $A$ and $B$ the degrees $2$ and $3\mspace{1mu}$, write $\Q[A,B]_n$ for the homogeneous part of degree~$n\mspace{1mu}$, and define $c_n\in\Q[A,B]_n$ by
\begin{equation}\label{eq:recurrence-c}
c_2\defis\frac A5\,,\qquad c_3\defis\frac B7\,,\qquad c_n\defis\frac{3}{(n-3)(2n+1)}\sum_{\substack{a+b=n\\a,b\geq2}}c_a\,c_b\qquad(n\geq4).
\end{equation}
Then $c_4=\frac{A^2}{75}\mspace{1mu}$, $c_5=\frac{3AB}{385}\mspace{1mu}$ and $c_6=\frac{2A^3}{4875}+\frac{B^2}{637}\mspace{1mu}$. By induction on~$n$ one checks that every monomial of degree~$n$ occurs in~$c_n$ with a positive coefficient and, writing $k=2n$ in~\eqref{eq:recurrence-G}, that
\begin{equation}\label{eq:specialization}
c_n\bigl(30\,G_4,\,70\,G_6\bigr)\=2(2n-1)\,G_{2n}\qquad(n\geq2)\,.
\end{equation}
Hence $c_n$ is, up to rescaling, the polynomial expressing $G_{2n}$ in terms of $G_4$ and~$G_6\mspace{1mu}$. The factors $30$ and~$70$ are chosen such that the generating series
\begin{equation}\label{eq:wpab}
\wp_{A,B}(z)\defis z^{-2}\+\sum_{n\geq2}c_n(A,B)\,z^{2n-2}
\end{equation}
satisfies the differential equation $\wp'^2=4(\wp^3-A\wp-B)$ of the Weierstrass $\wp$-function. 
In fact,  this differential equation is equivalent to~\eqref{eq:recurrence-c} by \cref{lem:ode}\ref{it:ode-unique}.

\begin{rk}\label{rk:elliptic} For comparison, if $\Lambda\subset\C$ is a lattice and $G_k(\Lambda)=\frac12\sum_{\omega\in\Lambda\smallsetminus\{0\}}\omega^{-k}\mspace{1mu}$, then $\wp_\Lambda=\wp_{A,B}$ for $A=30G_4(\Lambda)$ and $B=70G_6(\Lambda)\mspace{1mu}$, hence $c_n\bigl(30G_4(\Lambda),70G_6(\Lambda)\bigr)=2(2n-1)G_{2n}(\Lambda)\mspace{1mu}$. For $\Lambda=2\pi i(\Z\tau+\Z)$ we have $G_k(\Lambda)=G_k(\tau)\mspace{1mu}$, and this recovers~\eqref{eq:specialization}. These facts rely on the theory of elliptic functions, which we do not use.
\end{rk}

By~\eqref{eq:specialization}, the homomorphism $\Q[A,B]\to\Q[\![q]\!]$ given by $A\mapsto30G_4$ and $B\mapsto70G_6$ sends $c_a\,c_b\,c_c$ to a nonzero multiple of $G_{2a}G_{2b}G_{2c}\mspace{1mu}$, and $c_r\,c_{n-r}$ to a nonzero multiple of $G_{2r}G_{2n-2r}\mspace{1mu}$. Hence \cref{thm:triple} follows from the following statement about the~$c_n\mspace{1mu}$.

\begin{thm}\label{thm:polynomial}
For every $n\geq4\mspace{1mu}$, we have
\[
\Q[A,B]_n\=\spanQ\{c_a\,c_{n-a}\::\:2\leq a\leq n-2\}\,.
\]
\end{thm}

Conversely, \cref{thm:polynomial} follows from~\eqref{eq:Mk}, because $G_4$ and $G_6$ are algebraically independent.

\subsection*{The idea of the proof}
Let $e_1,e_2,e_3$ be the roots of $u^3-Au-B\mspace{1mu}$, and let $p_n\defis e_1^n+e_2^n+e_3^n\in\Z[A,B]$ be their power sums, so that $e_1e_2+e_1e_3+e_2e_3=-A$, $e_1e_2e_3=B$ and
\[
\sum_{n\geq0}p_n\,t^n\=\sum_{i=1}^3\frac{1}{1-e_it}\=\frac{3-At^2}{1-At^2-Bt^3}\,.
\]
Denote by $\Zt$ the ring of rational numbers with odd denominator. A set of polynomials in $\Zt[A,B]_n$ spans $\Q[A,B]_n$ as soon as its reductions modulo~$2$ span $\F[A,B]_n\mspace{1mu}$. The main observation is that the~$c_n\mspace{1mu}$, suitably normalized, are congruent to the~$p_n$ modulo~$2$, so that \cref{thm:polynomial} reduces to a statement about products of power sums modulo~$2\mspace{1mu}$. By \cref{prop:F2}, these span exactly the polynomials in $\F[A,B]_n$ that are divisible by~$B\mspace{1mu}$. For odd~$n$ this is all of $\F[A,B]_n\mspace{1mu}$, and for even~$n$ we require a further argument modulo~$4\mspace{1mu}$.

For $n\geq2$ let the normalized coefficient~$\tilde c_n$ be the multiple of~$c_n$ in which the coefficient of $A^{n/2}$ equals~$2$ if $n$ is even, and the coefficient of $A^{(n-3)/2}B$ equals~$1$ if $n$ is odd. This is well defined, since every monomial occurs in~$c_n$ with a positive coefficient.
\begin{table}[ht]
\centering
\renewcommand{\arraystretch}{1.3}
\begin{tabular*}{\textwidth}{@{\extracolsep{\fill}}c|ccccccc}
$n$ & $2$ & $3$ & $4$ & $5$ & $6$ & $7$ & $8$ \\ \hline
$c_n$ & $\frac15A$ & $\frac17B$ & $\frac1{75}A^2$ & $\frac3{385}AB$ & $\frac2{4875}A^3\+\frac1{637}B^2$ & $\frac2{5775}A^2B$ & $\frac1{82875}A^4\+\frac{12}{119119}AB^2$ \\
$\tilde c_n$ & $2A$ & $B$ & $2A^2$ & $AB$ & $2A^3\+\frac{375}{49}B^2$ & $A^2B$ & $2A^4\+\frac{9000}{539}AB^2$ \\
$p_n$ & $2A$ & $3B$ & $2A^2$ & $5AB$ & $2A^3\+3B^2$ & $7A^2B$ & $2A^4\+8AB^2$
\end{tabular*}
\\[2ex]
\begin{tabular*}{\textwidth}{@{\extracolsep{\fill}}c|ccc}
$n$ & $9$ & $10$ & $11$ \\ \hline
$c_n$ & $\frac{29}{2377375}A^3B\+\frac1{84721}B^3$ & $\frac2{6215625}A^5\+\frac{97}{19654635}A^2B^2$ & $\frac{778}{1991900625}A^4B\+\frac{246}{260275015}AB^3$ \\
$\tilde c_n$ & $A^3B\+\frac{1375}{1421}B^3$ & $2A^5\+\frac{181875}{5929}A^2B^2$ & $A^4B\+\frac{46125}{19061}AB^3$ \\
$p_n$ & $9A^3B\+3B^3$ & $2A^5\+15A^2B^2$ & $11A^4B\+11AB^3$
\end{tabular*}
\medskip
\caption{The coefficients $c_n$ of~\eqref{eq:recurrence-c}, their normalizations~$\tilde c_n$ and the power sums~$p_n\mspace{1mu}$.}\label{tab:Phi}
\vspace{-\baselineskip}
\end{table}

 In \cref{tab:Phi} the~$\tilde c_n$ have odd denominators, which holds in general:
\begin{prop}\label{prop:mod2}
For $n\geq2$ we have $\tilde c_n\in\Zt[A,B]$ and
\[
\tilde c_n\;\equiv\;p_n\pmod{2\Zt[A,B]}\,.
\]
\end{prop}

The main work of this paper lies in proving \cref{prop:mod2}. Its proof uses an operator~$\cT$ on $\Q[A,B]$ which has $c_n$ as an eigenvector with the odd eigenvalue $1+2^{2n-1}\mspace{1mu}$, and which, on a suitable $\Zt$-submodule of $\Q[A,B]_n\mspace{1mu}$, is congruent modulo~$2$ to a map of rank one. The operator~$\cT$ is defined by explicit substitutions in $A$ and~$B\mspace{1mu}$, and it has a classical meaning: for $A=30G_4(\Lambda)$ and $B=70G_6(\Lambda)$ as in \cref{rk:elliptic}, it is the Hecke operator~$T_2\mspace{1mu}$, which sums over the three lattices containing~$\Lambda$ with index~$2\mspace{1mu}$, i.e., over the three $2$-isogenies $\C/\Lambda\to\C/\Lambda'$ of the elliptic curve~$\C/\Lambda$ (see \cref{rk:isogeny}).

In \cref{sec:mod2} we prove the differential equation of~$\wp_{A,B}\mspace{1mu}$, use it to construct~$\cT\mspace{1mu}$, and prove \cref{prop:mod2}. In \cref{sec:proof} we study products of power sums modulo~$2$ and deduce \cref{thm:polynomial}.

\subsection*{Acknowledgments}
We thank Don Zagier for interesting suggestions and discussions. 

This project was partially supported by JSPS KAKENHI Grant 26K22254.
The second author's research is part of the project \emph{Modular forms and multiple zeta values in mathematical physics} with file number VI.Veni.242.073 of the research programme Veni, which is financed by the Dutch Research Council (NWO) under the grant \href{https://doi.org/10.61686/OJGCY95356}{doi.org/10.61686/OJGCY95356}.

\section{The normalized coefficients modulo 2}\label{sec:mod2}
In this section we prove \cref{prop:mod2}. As in the introduction, $A$ and $B$ are independent variables of degrees $2$ and~$3\mspace{1mu}$, and $c_n\in\Q[A,B]_n$ is defined by~\eqref{eq:recurrence-c}. For $\alpha,\beta$ in a field containing~$\Q(A,B)\mspace{1mu}$, let $Q_{\alpha,\beta}(u)\defis u^3-\alpha u-\beta\mspace{1mu}$, and write $Q=Q_{A,B}\mspace{1mu}$. Let $K$ be the splitting field of~$Q$ over $\Q(A,B)\mspace{1mu}$, and for $\alpha,\beta\in K$ let $\wp_{\alpha,\beta}\in K(\!(z)\!)$ be the series obtained by substituting $\alpha,\beta$ into the~$c_n\mspace{1mu}$ in \eqref{eq:wpab}. These series are characterized by the following differential equation.

\begin{lem}\label{lem:ode}
Let $\alpha,\beta\in K\mspace{1mu}$.
\begin{enumerate}[label=\upshape(\roman*)]
\item\label{it:ode-unique} For $f\in z^{-2}+K[\![z^2]\!]$ we have $f'^2=4\,Q_{\alpha,\beta}(f)$ if and only if $f=\wp_{\alpha,\beta}\mspace{1mu}$.
\item\label{it:ode-transform} Let $\varphi\in K(u)$ be such that $\varphi(u)-u\in K[\![u^{-1}]\!]\mspace{1mu}$, and
\[
\varphi'(u)^2\=\frac{Q_{\alpha,\beta}(\varphi(u))}{Q(u)}\,.
\]
Then $\varphi(\wp_{A,B})=\wp_{\alpha,\beta}\mspace{1mu}$.
\end{enumerate}
\end{lem}
\begin{proof}
Write $f=z^{-2}+\sum_{n\geq1}f_n\,z^{2n-2}$ for $f_n\in K\mspace{1mu}$. We first consider $f''-6f^2+2\alpha$. The coefficient of $z^{-2}$ in this expression is $-12f_1\mspace{1mu}$. If $f_1=0\mspace{1mu}$, the coefficient of $z^{2n-4}$ is $-10f_2+2\alpha$ for $n=2\mspace{1mu}$, zero for $n=3\mspace{1mu}$, and for $n\geq4$ it equals
\[
(2n-6)(2n+1)\,f_n\meno6\sum_{\substack{a+b=n\\a,b\geq2}}f_af_b\,.
\]
Hence $f''=6f^2-2\alpha$ if and only if $f_1=0\mspace{1mu}$, $f_2=\alpha/5$ and the $f_n$ with $n\geq4$ satisfy the recurrence~\eqref{eq:recurrence-c}, with $f_3$ arbitrary.

This second-order equation is related to the equation in~\ref{it:ode-unique} as follows: the derivative of $f'^2-4\,Q_{\alpha,\beta}(f)$ is $2f'(f''-6f^2+2\alpha)\mspace{1mu}$, and if $f_1=0\mspace{1mu}$, its constant term is $4\beta-28f_3\mspace{1mu}$. If $f'^2=4\,Q_{\alpha,\beta}(f)\mspace{1mu}$, then, as $f'$ is invertible in $K(\!(z)\!)\mspace{1mu}$, we get $f''=6f^2-2\alpha\mspace{1mu}$, and the constant term gives $f_3=\beta/7\mspace{1mu}$. Hence $f=\wp_{\alpha,\beta}\mspace{1mu}$. Conversely, $\wp_{\alpha,\beta}$ satisfies $f''=6f^2-2\alpha\mspace{1mu}$, so $\wp_{\alpha,\beta}'^2-4\,Q_{\alpha,\beta}(\wp_{\alpha,\beta})$ is constant, and its constant term is $4\beta-28c_3(\alpha,\beta)=0\mspace{1mu}$. This proves~\ref{it:ode-unique}.

For~\ref{it:ode-transform}, expand $\varphi$ in $K(\!(u^{-1})\!)$ and substitute $u=\wp_{A,B}\mspace{1mu}$. As $\wp_{A,B}^{-1}\in z^2K[\![z^2]\!]\mspace{1mu}$, this substitution is a well-defined homomorphism $K(\!(u^{-1})\!)\to K(\!(z)\!)$. By assumption on $\varphi$, we get $\varphi(\wp_{A,B})\in z^{-2}+K[\![z^2]\!]\mspace{1mu}$. By the chain rule and~\ref{it:ode-unique} for~$\wp_{A,B}\mspace{1mu}$,
\[
\bigl(\varphi(\wp_{A,B})'\bigr)^2\=\varphi'(\wp_{A,B})^2\cdot4\,Q(\wp_{A,B})\=4\,Q_{\alpha,\beta}\bigl(\varphi(\wp_{A,B})\bigr),
\]
and~\ref{it:ode-unique} for~$\wp_{\alpha,\beta}$ gives the claim.
\end{proof}

For a root~$e$ of~$Q\mspace{1mu}$, let
\[
 A_e\defis 15e^2-4A\,,\qquad B_e\defis 2e\,(11e^2-4A)
\]
and for $P\in \Q[A,B]$, define
\begin{equation}\label{eq:defT}
\cT(P)\defis\frac12\sum_eP(A_e,B_e),
\end{equation}
where the sum runs over the three roots $e$ of $Q$.
\begin{lem}\label{lem:tp}
Equation~\eqref{eq:defT} defines a linear map~$\cT: \Q[A,B]_n\to \Q[A,B]_n$ satisfying
\begin{equation*}
\cT(c_n)\=(1+2^{2n-1})\,c_n\qquad (n\geq 2).
\end{equation*}
\end{lem}
\begin{proof}
By Newton's identities the~$p_k$ lie in $\Z[A,B]\mspace{1mu}$, so that $\sum_e P(A_e,B_e)\in\Q[A,B]$ for every $P\in\Q[A,B]\mspace{1mu}$. As $A_e$ and~$B_e$ are homogeneous of degrees~$2$ and~$3$, $\cT$ preserves $\Q[A,B]_n\mspace{1mu}$.

We apply \cref{lem:ode}\ref{it:ode-transform} to the rational functions
\[
h_e(u)\defis u\+\frac{Q'(e)}{u-e}\,,\qquad \psi(u)\defis Q'(u)^2/Q(u)-8u\mspace{1mu},
\]
where $e$ is again a root of~$Q\mspace{1mu}$.
One checks directly, using $Q(e)=0\mspace{1mu}$, that
\[
h_e'^2\=\frac{Q_{A_e,B_e}(h_e)}{Q}\,,\qquad \psi'^2\=\frac{Q_{16A,64B}(\psi)}{Q}\,.
\]
Since $h_e(u)-u$ and $\psi(u)-u$ lie in $u^{-1}K[\![u^{-1}]\!]\mspace{1mu}$, \cref{lem:ode}\ref{it:ode-transform} gives
\[
\wp_{A_e,B_e}\=h_e(\wp_{A,B}),\qquad \wp_{16A,64B}\=\psi(\wp_{A,B}).
\]
The partial fraction decomposition $Q'(u)^2/Q(u)=9u+\sum_eQ'(e)/(u-e)\mspace{1mu}$, where $e$ runs over the three roots of~$Q\mspace{1mu}$, shows that $\sum_eh_e(u)=\psi(u)+2u\mspace{1mu}$. Hence 
\[\sum_e\wp_{A_e,B_e}\=2\wp_{A,B}+\wp_{16A,64B}\,.\]
As $c_n$ is homogeneous of degree~$n\mspace{1mu}$, we have $c_n(16A,64B)=4^nc_n\mspace{1mu}$, and comparing the coefficients of~$z^{2n-2}$ gives $\sum_ec_n(A_e,B_e)=(2+4^n)\,c_n\mspace{1mu}$.
\end{proof}

\begin{rk}\label{rk:isogeny}
Let $A=30G_4(\Lambda)$ and $B=70G_6(\Lambda)$ for a lattice $\Lambda\subset\C\mspace{1mu}$, so that $\wp_{A,B}=\wp_\Lambda\mspace{1mu}$. The roots of~$Q$ are the values $e=\wp_\Lambda(\omega/2)$ with $\omega\in\Lambda\setminus2\Lambda\mspace{1mu}$, and $h_e(\wp_\Lambda(z))=\wp_\Lambda(z)+\wp_\Lambda(z+\omega/2)-e$ is the $\wp$-function of the lattice $\Lambda_e=\Lambda+\Z\,\omega/2\mspace{1mu}$, one of the three lattices containing~$\Lambda$ with index~$2\mspace{1mu}$.  So $A_e=30G_4(\Lambda_e)$ and $B_e=70G_6(\Lambda_e)\mspace{1mu}$. These are V\'elu's formulas~\cite{V} for $2$-isogenies. Similarly, $\psi(\wp_\Lambda(z))=4\wp_\Lambda(2z)=\wp_{\Lambda/2}(z)$ is the duplication formula. In this language $\cT$ is the Hecke operator~$T_2\mspace{1mu}$, and \Cref{lem:tp} follows by counting: each of the three lattices~$\Lambda_e$ consists of~$\Lambda$ and one of the three nonzero cosets of~$\Lambda$ in~$\frac12\Lambda\mspace{1mu}$, so every element of $\frac12\Lambda\setminus\Lambda$ lies in exactly one~$\Lambda_e$ and every element of~$\Lambda$ in all three. Hence
\[
\sum_eG_{2n}(\Lambda_e)\=G_{2n}\bigl(\tfrac12\Lambda\bigr)\+2\,G_{2n}(\Lambda)\=(4^n+2)\,G_{2n}(\Lambda)\,.
\]
\end{rk}

We now reduce modulo~$2\mspace{1mu}$. Write $R=\Zt[A,B]$ and $R_n$ for its homogeneous part of degree~$n\mspace{1mu}$. As $Q$ is monic, $\sum_eF(e)\in R$ for $F\in R[u]\mspace{1mu}$, and $\sum_eF(e)\in2^jR$ if $F\in2^jR[u]\mspace{1mu}$. The proof of \cref{prop:mod2} rests on two facts: $1+2^{2n-1}$ is odd, and
\[
A_e\;\equiv\;e^2\pmod 2\,,\qquad B_e\;\equiv\;2e^3\pmod 4\,.
\]

\begin{proof}[Proof of \cref{prop:mod2}]

First we choose as a basis of $\Q[A,B]_n\mspace{1mu}$ the monomials $A^aB^b$ of degree~$n$ with $a\geq 0$ and $b\geq2\mspace{1mu}$, together with $v_0=A^{(n-3)/2}B$ if $n$ is odd and $v_0=2A^{n/2}$ if $n$ is even. Let $L_n$ be its $\Zt$-span. Then, $L_n=R_n$ if $n$ is odd, and $L_n$ consists of the polynomials in $R_n$ whose coefficient of~$A^{n/2}$ is even if $n$ is even. By definition, for $\tilde c_n$ the coefficient of $v_0$ in this basis equals~$1$.

We claim that $\cT(v)\in2L_n$ for every basis vector $v\neq v_0\mspace{1mu}$, and that $\cT(v_0)\in L_n$ with $\cT(v_0)\equiv p_n\pmod{2R_n}\mspace{1mu}$. For $v=A^aB^b$ with $a\geq 0$ and $b\geq2\mspace{1mu}$, the sum $\sum_eA_e^aB_e^b$ lies in~$2^bR\mspace{1mu}$, so $\cT(v)\in2^{b-1}R_n\subset2R_n\mspace{1mu}$. This proves the first claim for odd~$n\mspace{1mu}$. For $n=2m$ even we also need the coefficient of~$A^m$ to be divisible by~$4\mspace{1mu}$. We read it off at $B=0\mspace{1mu}$, where $Q$ specializes to $u^3-Au$ with roots $0$ and~$\pm\sqrt A\mspace{1mu}$, giving $(A_e,B_e)$ equal to $(-4A,0)$ and $(11A,\pm14A\sqrt A)$, respectively. As $b\geq 2$ is even, this gives
\[
[A^m]\,\cT(A^aB^b)\=11^a\,14^b\;\in\;4\Z\,.
\]
For~$v_0\mspace{1mu}$, we have $\cT(v_0)=\sum_eA_e^{(n-3)/2}e(11e^2-4A)$ if $n$ is odd and $\cT(v_0)=\sum_eA_e^{n/2}$ if $n$ is even. In both cases the summand is congruent to~$e^n$ modulo~$2\mspace{1mu}$, so that $\cT(v_0)\equiv p_n\pmod{2R_n}\mspace{1mu}$. For odd~$n$ this also shows that $\cT(v_0)\in L_n\mspace{1mu}$. For $n=2m$ even, evaluating at $B=0$ as before gives $[A^m]\cT(v_0)=(-4)^m+2\cdot11^m\in2\Z\mspace{1mu}$, so again $\cT(v_0)\in L_n\mspace{1mu}$.

By the claim, $\cT$ preserves~$L_n\mspace{1mu}$ and $\cT(P)\equiv x_0\,\cT(v_0)\pmod{2L_n}$ for every $P\in L_n$ with $v_0$-coordinate~$x_0\mspace{1mu}$. Now choose a rational multiple~$P$ of~$c_n$ with $P\in L_n\setminus2L_n\mspace{1mu}$. As $1+2^{2n-1}$ is odd, \Cref{lem:tp} gives
\[
P\;\equiv\;(1+2^{2n-1})\,P\=\cT(P)\;\equiv\;x_0\,\cT(v_0)\pmod{2L_n}\,,
\]
so that $x_0$ is not divisible by~$2$ in~$\Zt\mspace{1mu}$. Hence $\tilde c_n=P/x_0$ lies in $L_n\subset R_n\mspace{1mu}$, and the same congruence with~$\tilde c_n$ in place of~$P$ gives $\tilde c_n\equiv\cT(v_0)\equiv p_n\pmod{2R_n}\mspace{1mu}$.
\end{proof}

\section{Proof of \texorpdfstring{\cref{thm:polynomial}}{Theorem 3}}\label{sec:proof}
We first prove the analogue of \cref{thm:polynomial} for the power sums over~$\F\mspace{1mu}$, and then lift it using \cref{prop:mod2}.

\begin{prop}\label{prop:F2}
For every $n\geq3\mspace{1mu}$, the reductions of the power sums modulo~$2$ satisfy
\[
\spanF\{p_n,\,p_ap_{n-a}\::\:2\leq a\leq n-2\}\=B\,\F[A,B]_{n-3}\,.
\]
\end{prop}
\begin{proof}
Over~$\F$ we have $u^3-Au-B=(u-x)(u-y)(u-x-y)$ in $\F[x,y]\mspace{1mu}$, where
\[
A\=x^2+xy+y^2\,,\qquad B\=xy(x+y)\,.
\]
Since $x$ and~$y$ are algebraic over $\F(A,B)\mspace{1mu}$, the polynomials $A$ and~$B$ are algebraically independent. We may therefore view $\F[A,B]$ as a subring of $\F[x,y]\mspace{1mu}$, in which $p_m=x^m+y^m+(x+y)^m\mspace{1mu}$. For $f\in\F[x,y]$ let
\[
\tau(f)(x,y)\defis f(x,y)\+f(y,x+y)\+f(x+y,x)
\]
be the sum over the cyclic permutations of the three roots. As $A$ and~$B$ are invariant under all permutations of the roots and $3=1$ in~$\F\mspace{1mu}$, we have $\tau(f)=f$ for $f\in\F[A,B]\mspace{1mu}$. On the other hand, for $s_i:=x^iy^{n-i}+x^{n-i}y^i$ the six terms of $\tau(s_i)$ are the products $e^i\tilde e^{\,n-i}$ over all ordered pairs of distinct roots $e,\tilde e\mspace{1mu}$, so that
\[
\tau(s_i)\=p_i\,p_{n-i}\meno p_n\=p_i\,p_{n-i}\+p_n\,.
\]
This is zero for $i=0$ (as $p_0=1$ in~$\F$) and equals~$p_n$ for $i=1$ (as $p_1=0$).

An element $f\in B\,\F[A,B]_{n-3}$ is symmetric in $x$ and~$y$ and vanishes at $x=y\mspace{1mu}$. It is therefore a linear combination of the~$s_i$ with $0\leq i<n/2\mspace{1mu}$, since the middle monomial $x^{n/2}y^{n/2}$ does not vanish at $x=y\mspace{1mu}$. Applying~$\tau$ expresses $f=\tau(f)$ as a linear combination of $p_n$ and the $p_ip_{n-i}$ with $2\leq i\leq n-2\mspace{1mu}$. Conversely, for $m\geq1$ the power sum~$p_m$ vanishes at $x=y\mspace{1mu}$, where the roots are $x,x,0\mspace{1mu}$. An element $f\in\F[A,B]$ vanishing at $x=y$ satisfies $f(x^2,0)=0\mspace{1mu}$, as $A=x^2$ and $B=0$ there, and hence is divisible by~$B\mspace{1mu}$. So $p_n$ and the $p_ap_{n-a}$ lie in $B\,\F[A,B]_{n-3}\mspace{1mu}$.
\end{proof}
For example, for $n=6$ we have $p_2p_4\equiv0$ and $p_6\equiv p_3^2\equiv B^2\mspace{1mu}$. So by \cref{prop:mod2}, the reductions modulo~$2$ of $\tilde c_6\mspace{1mu}$, $\tilde c_2\tilde c_4$ and~$\tilde c_3^2$ span
\[
\spanF\{p_6,\,p_2p_4,\,p_3^2\}\=\F B^2\;\not\ni\;A^3\,,
\]
and we cannot conclude directly that these three polynomials span $\Q[A,B]_6\mspace{1mu}$. The proof below overcomes this by working modulo~$4\mspace{1mu}$: $\tilde c_6$ and~$\tilde c_3^2$ are congruent modulo~$2$ but not modulo~$4\mspace{1mu}$, and
\[
\tfrac12\bigl(\tilde c_6-\tilde c_3^2\bigr)\=A^3+\tfrac{163}{49}B^2
\]
has odd denominators and reduces modulo~$2$ to $A^3+B^2\mspace{1mu}$.

\begin{proof}[Proof of \cref{thm:polynomial}]
By~\eqref{eq:recurrence-c}, $c_n$ is a linear combination of the $c_ac_{n-a}\mspace{1mu}$, so the right-hand side of \cref{thm:polynomial} is the $\Q$-span~$V_n$ of~$\tilde c_n$ and the products $\tilde c_a\tilde c_{n-a}$ with $2\leq a\leq n-2\mspace{1mu}$. Let $V_n^{(2)}\subset R_n$ be the $\Zt$-span of these polynomials, and, if $n$ is even, of
\[
\Psi_n\defis \frac{\tilde c_n-\tilde c_{n/2}^2}{2}\,.
\]
It lies in~$V_n\mspace{1mu}$, as $\tilde c_{n/2}^2$ is one of the products. It also lies in~$R_n\mspace{1mu}$, since $\tilde c_n\equiv p_n\equiv p_{n/2}^2\equiv\tilde c_{n/2}^2$ modulo~$2$ by \cref{prop:mod2}. Its coefficient of~$A^{n/2}$ is $\pm1\mspace{1mu}$: we have $\tilde c_n(A,0)=2A^{n/2}\mspace{1mu}$, whereas $\tilde c_{n/2}(A,0)$ equals $2A^{n/4}$ if $n/2$ is even and~$0$ if $n/2$ is odd. 

By \cref{prop:mod2,prop:F2}, the reduction of~$V_n^{(2)}$ modulo~$2$ contains $B\,\F[A,B]_{n-3}\mspace{1mu}$. For odd~$n$ this is all of $\F[A,B]_n\mspace{1mu}$, since every monomial of degree~$n$ contains~$B\mspace{1mu}$. For even~$n$ the reduction of~$\Psi_n$ supplies the missing monomial~$A^{n/2}\mspace{1mu}$. In both cases $V_n^{(2)}$ reduces onto $\F[A,B]_n\mspace{1mu}$. Hence some maximal minor of the matrix of coefficients of its generators is a unit in~$\Zt\mspace{1mu}$, and in particular nonzero. Therefore $V_n=\Q[A,B]_n\mspace{1mu}$.
\end{proof}

\begin{rk}\label{rk:periods}
In a previous approach, we aimed to generalize the period polynomial relations satisfied by generating series of products of two Eisenstein series to the generating series of triple products. Modulo the~$G_k\mspace{1mu}$, the generating series $\sum_{a+b=k}G_aG_b\,x^{a-1}y^{b-1}\mspace{1mu}$, suitably completed by terms involving~$G_2\mspace{1mu}$, is a period polynomial with coefficients in~$M_k\mspace{1mu}$, i.e., it satisfies (see \cite{GKZ}*{Theorems~4 and~5})
\[
f(x,y)\+f(-y,x)\=0\,,\qquad f(x,y)\+f(x-y,x)\+f(y,y-x)\=0\,.
\]
For the generating series
\[
F(x,y,z)\=\sum_{a,b,c}G_a\,G_b\,G_c\;x^{a-1}y^{b-1}z^{c-1}
\]
of triple products, suitably completed by terms involving~$G_2\mspace{1mu}$, we hoped that, modulo the span of the~$G_k$ and the products of two Eisenstein series, it satisfies the ``three-variable period polynomial relations''
\begin{equation}\label{eq:3rel}
f(x,y,z)\=f(y,z,x)\,,\qquad f(x,y,z)\+f(-y,x-y,z)\+f(y-x,-x,z)\=0\,.
\end{equation}
The first relation holds since $F$ is symmetric. Based on numerical evidence, we conjectured that every $f\in\Q[x,y,z]$ satisfying both relations lies in $\Q[x,y]+\Q[y,z]+\Q[z,x]\mspace{1mu}$. This was recently proved by Ilin and Ono, with a proof found in collaboration with the AI system AxiomProver. It follows that every solution is of the form $g(x,y)+g(y,z)+g(z,x)$ for a period polynomial~$g\mspace{1mu}$, and hence that there are no nonzero solutions of odd degree \cite{IO}*{Theorem~1.1 and Corollary~1.3}. Since $F$ is homogeneous of odd degree in each weight, this would give $F\equiv0\mspace{1mu}$, i.e., \cref{thm:triple}. However, it seems that showing that~$F$ satisfies \eqref{eq:3rel} is as hard as \cref{thm:triple} itself.
\end{rk}

\subsection*{Use of AI tools} This paper is the outcome of a long-running project on period polynomials and products of Eisenstein series, which contained several approaches and partial results to an elementary proof of~\eqref{eq:Mk}. AI models helped fill the missing key step for one of these approaches. We then wrote this paper ourselves, using LLMs for language editing only.

\end{document}